\documentclass[11pt]{article} 
 
\usepackage{epsfig}
\usepackage{amssymb}
\usepackage[]{times} 
\newcommand{\be}{\begin{equation}}
\newcommand{\ee}{\end{equation}} 
\newcommand{\bea}{\begin{eqnarray}}
\newcommand{\eea}{\end{eqnarray}}

\def\1#1{^{(#1)}}
\def\la{\langle} 
\def\ra{\rangle}
\begin{document} 
\title{On Frobenius numbers for Symmetric\\ (not Complete Intersection) 
\\Semigroups Generated by Four Elements} \author{Leonid G. Fel\\ \\ 
Department of Civil Engineering, Technion, Haifa 32000, Israel} 
\date{} 
\maketitle
\def\be{\begin{equation}}
\def\ee{\end{equation}} 
\def\bea{\begin{eqnarray}}
\def\eea{\end{eqnarray}}
\def\p{\prime} 
\vspace{-1cm} 
\begin{abstract} 
We derive the lower bound for Frobenius number of symmetric (not complete 
intersection) semigroups generated by four elements.
\end{abstract} 
\section{Introduction}\label{l1}
Bresinsky \cite{be752} has shown that symmetric numerical semigroups ${\sf S}_4
=\la d_1,d_2,d_3,d_4\ra$ which are not complete intersection (CI), have the 
first Betti number $\beta_1=5$, i.e., its Hilbert series $H\left({\sf S}_4;z
\right)$ reads, 
\bea 
H\left({\sf S}_4;z\right)=\frac{1-\sum_{j=1}^5z^{a_j}+\sum_{j=1}^5z^{c-a_j}-z^
{c}}{\prod_{i=1}^4\left(1-z^{d_i}\right)},\quad a_j,c\in{\mathbb Z}_+,\label{x1}
\eea
where $\gcd(d_1,d_2,d_3,d_4)=1$ and $d_i>4$; the latter is necessary since the
numerical semigroup $\la m,d_2,\ldots,d_m\ra$ is never symmetric \cite{f11}. 

Denote the k-th power symmetric polynomials in five variables $I_k(a_1,\ldots,
a_5)=\sum_{j=1}^5a_j^k$ for short by $I_k$ and write two polynomial identities 
derived in \cite{fel09}, formula (6.11), 
\bea
8I_3-6I_2I_1+I_1^3=24\pi_4\;,\quad 2c=I_1\;,\quad\pi_4=d_1d_2d_3d_4\;.\label{x2}
\eea 
In the present paper we show that relations (\ref{x2}) allow to find a strong 
lower bound to the Frobenius number $F_{\widetilde{CI}}({\sf S}_4)$ of symmetric
(not CI) semigroup ${\sf S}_4$.

Consider elementary symmetric polynomials $J_k=J_k(a_1,\ldots,a_5)$ in five 
variables
\bea
J_0=1,\quad J_1=\!\sum_{j=1}^5a_j,\quad J_2=\!\sum_{i\geq j=1}^5a_ia_j,\quad 
J_3=\!\sum_{i\geq j=1}^5\frac{J_5}{a_ia_j},\quad J_4=\!\sum_{j=1}^5\frac{J_5}
{a_j},\quad J_5=\prod_{j=1}^5a_j,\label{x3}
\eea
which are related to each other by inequalities \cite{har59},
\bea
\left(\frac{J_r}{{5\choose r}}\right)^{\frac{1}{r}}\geq\left(\frac{J_{r+1}}
{{5\choose r+1}}\right)^{\frac{1}{r+1}},\quad\mbox{i.e.}\quad\frac{J_1}{5}\geq
\sqrt{\frac{J_2}{10}}\geq\sqrt[3]{\frac{J_3}{10}}\geq\sqrt[4]{\frac{J_4}{5}}\geq
\sqrt[5]{J_5}.\label{x4}
\eea
Make use of the Newton recursion identities for symmetric polynomials $I_r$ and 
$J_r$
\bea
&&mJ_m=\sum_{k=1}^m(-1)^{k-1}I_kJ_{m-k},\quad\mbox{i.e.}\label{x5}\\
&&I_1=J_1,\quad I_2=J_1I_1-2J_2,\quad I_3=J_1I_2-J_2I_1+3J_3,\quad\mbox{etc},
\nonumber
\eea
and write explicit expressions for the three first of them, $m=1,2,3$,
\bea
I_1=J_1,\quad I_2=J_1^2-2J_2,\quad I_3=J_1^3-3J_2J_1+3J_3.\label{x6}
\eea
Substitute (\ref{x6}) into (\ref{x2}) and rewrite the latter identities in 
$J_r$ polynomials,
\bea
J_1^3-4J_2J_1+8J_3=8\pi_4,\quad J_1=2c\;.\label{x7}
\eea
Combine two identities: (\ref{x7}) and $J_1\geq5\sqrt[3]{J_3/10}$ in (\ref{x4}),
and get
\bea
cJ_2+\pi_4=c^3+J_3\leq\frac{41}{25}c^3.\label{x8}
\eea
Thus we arrive at the first inequality in (\ref{x4}) and inequality (\ref{x8}),
\bea
a)\;\;J_2\leq\frac{8}{5}c^2,\hspace{.6cm}b)\;\;J_2\leq\frac{41}{25}c^2-\frac{
\pi_4}{c}.\label{x9}
\eea
Inequality (\ref{x9}a) holds always while inequality (\ref{x9}b) is valid not 
for every $c$. In order to make both inequalities consistent we have to find 
such range of $c$ where both inequalities (\ref{x9}) are satisfied for any $c$ 
within the range. A simple analysis shows that it can be provided iff 
\bea
c\geq c_*,\quad\mbox{where}\quad\frac{41}{25}c_*^2-\frac{\pi_4}{c_*}=\frac{8}{5}
c_*^2,\quad\mbox{i.e.}\quad c\geq\sqrt[3]{25\pi_4},\quad\sqrt[3]{25}\simeq 
2.924.\label{x10}
\eea
The last expression leads immediately to the lower bound $F_{\widetilde{CI}}
({\sf S}_4)$ of the Frobenius number of symmetric (not CI) semigroup generated 
by four elements,
\bea
F({\sf S}_4)\geq F_{\widetilde{CI}}({\sf S}_4)=\sqrt[3]{25\pi_4}-\sigma_4,\quad
\sigma_n=\sum_{i=1}^nd_i.\label{x11}
\eea
If the values of generators are not large $d_j\sim 10$ then the bound 
(\ref{x11}) is very close to the exact value of $F({\sf S}_4)$, e.g., for 
symmetric not CI numerical semigroups treated in \cite{de76,bfs14,fel09} we have
\bea
F(5,6,7,8)=9>8.76,\quad F(7,8,9,13)=19>17.715,\quad F(8,13,15,17)=35>34.198.
\nonumber
\eea
However, for large generators $d_j>100$ such proximity is broken. We present the
Hilbert series and Frobenius number of symmetric not CI numerical semigroup 
$\overline{S_4}=\la 151,154,157,158\ra$ found by author,
\bea
H\left(\overline{S_4};z\right)=\frac{1-z^{308}-z^{625}-z^{628}+z^{779}+z^{782}-
z^{3473}-z^{3476}+z^{3627}+z^{3630}+z^{3947}-z^{4255}}{\left(1-z^{151}\right)
\left(1-z^{154}\right)\left(1-z^{157}\right)\left(1-z^{158}\right)},\nonumber\\
F\left(\overline{S_4};z\right)=3635,\quad  F_{\widetilde{CI}}\left(\overline{S_4
};z\right)=1814.1\;.\hspace{5cm}\nonumber
\eea
This phenomenon is known for nonsymmetric numerical semigroups ${\sf S}_3$  
generated by three elements. The lower bound $F_{NS}({\sf S}_3)$ of Frobenius 
number for these semigroups was found in \cite{ro90,da94} and improved slightly 
in \cite{fel06}, $F_{NS}({\sf S}_3)=\sqrt{3}\sqrt{d_1d_2d_3+1}-\sigma_3$, e.g.,
\bea
F(\la 5,6,7\ra)=9>7.16,\quad F(\la 8,9,13\ra)=28>23,\quad F(\la 151,154,157\ra)
=11624\gg 2847.5\;.\nonumber
\eea
We finish by comparison the bound (\ref{x11}) with two other lower bounds of 
Frobenius numbers for symmetric CI \cite{fel09} $F_{CI}({\sf S}_4)$ and 
nonsymmetric \cite{kil00} $F_{NS}({\sf S}_4)$ semigroups generated by four 
elements,
\bea
F_{CI}({\sf S}_4)=3\sqrt[3]{\pi_4}-\sigma_4,\quad F_{NS}({\sf S}_4)=\sqrt[3]{6
\pi_4}-\sigma_4,\quad \sqrt[3]{6}\simeq 1.817.\nonumber
\eea
\section*{Acknowledgement}
The research was supported by the Kamea Fellowship.

\end{document}